\documentclass[12pt,a4paper]{article}

\usepackage{graphicx}
\usepackage{amsmath}
\usepackage{amssymb}
\usepackage{float}
\usepackage{url}
\usepackage{hyperref}

\newcommand{\smt}{\texttt{smt}}
\newcommand{\smcirc}{\texttt{smcirc}}
\newcommand{\smtoep}{\texttt{smtoep}}
\newcommand{\bm}[1]{\mathbf{#1}}
\newcommand{\ttt}[1]{\texttt{#1}}
\newcommand{\e}{{\mathrm{e}}}
\newcommand{\ii}{{\mathrm{i}}}

\newcommand{\cF}{{\mathcal{F}}}

\newcommand{\deltab}{\boldsymbol\delta}

\newtheorem{Algorithm}{Algorithm}
\newcounter{linea}[Algorithm]

\DeclareMathOperator{\diag}{diag}

\DeclareMathOperator{\FFT}{\texttt{fft}}
\DeclareMathOperator{\IFFT}{\texttt{ifft}}

\numberwithin{equation}{section}

\newfloat{falgo}{tbp}{loa}[section]
\floatname{falgo}{Algorithm}

\newcommand{\IF}{\textsf{if}\ }
\newcommand{\ELSE}{\textsf{else}\ }

\title{\smt: a Matlab structured matrices toolbox}
\author{Michela Redivo-Zaglia\thanks{%
Dipartimento di Matematica Pura e Applicata, via Trieste 63, 35121 Padova,
Italy (\texttt{Michela.RedivoZaglia@unipd.it}).} \footnotemark[3]
\and
Giuseppe Rodriguez\thanks{%
Dipartimento di Matematica e Informatica, Universit\`a di Cagliari, viale
Merello 92, 09123 Cagliari, Italy (\texttt{rodriguez@unica.it}).}
\footnote{This work was supported in part by MIUR, under the PRIN grant no.\
2006017542-003, and by the University of Padova Project no.\ CPDA089040.}}
\date{}

\begin{document}

\maketitle

\begin{abstract}
We introduce the \smt\ toolbox for Matlab.
It implements optimized storage and fast arithmetics for circulant and Toeplitz
matrices, and is intended to be transparent to the user and easily extensible.
It also provides a set of test matrices, computation of circulant
preconditioners, and two fast algorithms for Toeplitz linear systems.
\end{abstract}

\section{Introduction}\label{sec:intro}

Algebraic structures are present in many mathematical problems,
so they arise naturally in a large number of applications, like
medical imaging, remote sensing, geophysical prospection, image
deblurring, etc. Moreover, in many real-world computations,
the full exploitation of the structure of the problem is
essential
to be able to manage large dimensions and real time
processing.

In the last 25 years, a great effort has been made to study the properties of
algebraic structures and to develop algorithms capable of taking advantage
of
these structures in the solution of various matrix problems (solution of linear
systems, eigenvalues computation, etc.), as well as in matrix arithmetics, for
what concerns memory storage, speed of computation and stability.

In spite of the many important advances in this field, there is
not much software publicly available for structured matrices computation.
On the contrary, most of the \emph{fast} algorithms which have been proposed,
and whose properties have been studied theoretically, exist only under
the form of published papers or, in some occasion, unreleased (and often
unoptimized) research code.
This fact often force researchers to re-implement from scratch algorithms and
\emph{BLAS-like} routines, even for the most classical classes of structured
matrices.
Anyway, this is possible only for those with enough knowledge in Mathematics
and Computer Science, and totally rules out a large amount of potential users of
structured algorithms.

Matlab~\cite{matlab77} is a computational environment which is extremely
diffused among both applied mathematicians and engineers, in academic as well
as in industrial research.
It makes matrix computation sufficiently easy and immediate, and provides the
user with powerful scientific visualization tools.
At the moment, besides the standard unstructured (or \emph{full}) matrices, the
only matrix structure natively supported in Matlab is sparsity.
In \emph{sparse matrix storage} only nonzero elements are kept in memory,
together with their position inside the matrix.
Moreover, all operations between sparse matrices are redefined to reduce
execution time and memory consumption, while mixed computations return full or
sparse arrays depending on the type of operations involved.

Our idea is to extend Matlab with a computational framework devoted to structured
matrices, with the aim of making it easy to use and \emph{Matlab-like},
transparent for the user, highly optimized for what concerns storage and
complexity, and easily extensible.
We tried to follow closely the way Matlab treats sparse matrices, and for doing
this, we used the Matlab object-oriented classes.
Starting from version 5, in fact, it is possible to add new data types in
Matlab (\emph{classes}), to define \emph{methods} for \emph{classes}, i.e.,
functions to create and manipulate variables belonging to the new data
type, and to \emph{overload} (redefine) the arithmetic operators for each
new class.

At the moment, our toolbox supports two very common classes of structured
matrices, namely circulant and Toeplitz matrices.
In writing the software our aim was not only to furnish storage support,
full arithmetics and some additional methods for these structured matrices,
but also to create a framework easily extendible, in terms of functions and new
data types, and
to specify a pattern for future developments of the package.
So, a great effort was spent in the software engineering of the
toolbox.

Among the available Matlab software for structured matrices computation, we
mention the following Internet resources.
Various fast and superfast algorithms for structured matrices have been
developed by the MaSe-team (Matrices having Structure)~\cite{maseteam},
coordinated by Marc Van Barel at the Katholieke Universiteit of Leuven.
A Toolbox for Structured Matrix Decompositions~\cite{tsmd} has been included in
the SLICOT package~\cite{slicot}, developed under the NICONET (Numerics in
Control Network) European Project~\cite{niconet}.
The RestoreTools~\cite{restoretools} is an object oriented Matlab package for
image restoration which has been developed by James Nagy and his group at Emory
University.
The MOORe Tools~\cite{moore}, an object oriented toolbox for the solution of
discrete ill-posed problems derived from~\cite{han94}, provides some support
for certain classes of structured matrices, mainly Kronecker products,
circulant and block-circulant matrices.

Matlab implementations of various algorithms are also available in the
personal home pages of many researchers working in this field.
Many subroutines written in general purpose languages, like C or Fortran,
are also available.
It is worth mentioning that there are plans to add support for structured
matrices in LAPACK and ScaLAPACK; see~\cite[Section~4]{dd05}.

The plan of this paper is the following.
In Section~\ref{sec:toolbox}, we describe in detail our toolbox, called
\smt\ (Structured Matrix Toolbox), its capabilities, the new data types added
to Matlab and the functions for their treatment.
Section~\ref{sec:implementation} is devoted to some technical implementation
issues, while in Section~\ref{sec:conclusion} we describe possible future lines
of development of this software package.

\section{The toolbox}\label{sec:toolbox}

\begin{figure}[htbp]
\centerline{\includegraphics[width=.8\textwidth]{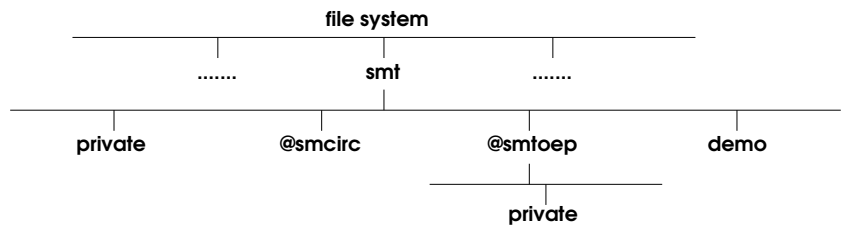}}
\caption{Directory tree of \smt}
\label{fig:dirtree}
\end{figure}

Once installed (see Section~\ref{sec:implementation}), the toolbox resides in
the directory tree sketched in Fig.~\ref{fig:dirtree}.
The main directory contains a set of general purpose functions, described in
detail in Section~\ref{sec:general}, and the following four subdirectories:
\begin{itemize}
\item \ttt{~@smcirc} and \ttt{@smtoep}, which contain the functions
to create and manipulate the objects of class \smcirc\ and
\smtoep, i.e., circulant and Toeplitz matrices;
\item \ttt{~private}, whose functions, discussed
in Section~\ref{sec:general}, are
accessible by the user only through the commands placed at the upper
directory level;
\item \ttt{~demo}, which hosts an interactive tutorial on the basic use of the
toolbox.
\end{itemize}

Let us briefly explain how Matlab deals with new data types.
When the user creates an object of class, say, \ttt{obj}, then the interpreter looks
for the function with the same name in a directory called \ttt{@obj}, located
in the search path.
Similarly, when an expression involves a variable of class \ttt{obj} or a
function is applied to it, the same directory is searched for an appropriate
operator or function defined for objects of this class.

%

Writing this software, we took great care in checking the validity of the input
parameters, in particular for what concerns dimensions and data types,
and in using an appropriate style for warnings and errors,
in order to guarantee
the \emph{Matlab-like} behaviour of the toolbox.
As this requires a long chain of conditional tests, the resulting functions are
often more complicated than expected (see for example the file \ttt{mtimes.m}
in the directory \ttt{@smcirc}), but this does not seem to have a significant
impact on execution time.

Full documentation for every function of the toolbox is accessible via the
Matlab \ttt{help} command and the code itself is extensively commented.
Manual pages can be obtained by the usual Matlab means,
i.e.,
\begin{quote}
\footnotesize
\begin{tabbing}
XXXXXXXXXXXXXXXXXXXX \= \kill
\ttt{help} $<$func\_name$>$ \> for the functions in the main directory, \\
\ttt{help} $<$class$>$/$<$func\_name$>$ \> where $<$class$>$ is either
    \ttt{@smcirc} or \ttt{@smtoep}, \\
\ttt{help private}/$<$func\_name$>$ \> for the functions in the \ttt{private}
    subdirectory.
\end{tabbing}
\end{quote}
Notice that $<$func\_name$>$ may be \ttt{Contents} (except in conjunction with
\ttt{private}), in which case a description of the entire directory content is
displayed.
For example, the command
\begin{quote}
\footnotesize
\begin{verbatim}
help @smtoep/Contents
\end{verbatim}
\end{quote}
displays the list of all the functions, operators and methods for 
\smtoep\ objects (i.e., Toeplitz matrices), while
\begin{quote}
\footnotesize
\begin{verbatim}
help @smtoep/mtimes
\end{verbatim}
\end{quote}
gives information about the matrix product operator for Toeplitz matrices.

We remark that the toolbox supports both real and complex structured matrices,
since complex numbers are natively implemented in Matlab.
It is also possible to manage sparse \smcirc\ and \smtoep\ objects.

\subsection{Circulant matrices}\label{sec:smcirc}

A circulant matrix~\cite{dav79} of order $n$ is a matrix $C$ whose elements
satisfy the relations
$$
\begin{aligned}
C_{ij} & = c_{i-j}, \quad & i,j&=1,\ldots,n, \\
c_{k-n} &= c_k, \quad & k&=1,\ldots,n-1,
\end{aligned}
$$
e.g., for $n=4$,
$$
C = \begin{bmatrix}
c_0 & c_3 & c_2 & c_1 \\
c_1 & c_0 & c_3 & c_2 \\
c_2 & c_1 & c_0 & c_3 \\
c_3 & c_2 & c_1 & c_0
\end{bmatrix}.
$$

The main property of a circulant matrix is that it is diagonalized by the
normalized Fourier matrix $\cF_\eta$, defined by
$$
(\cF_\eta)_{ij} = \frac{1}{\sqrt{n}} \eta^{ij},
$$
where $\eta$ is any primitive complex $n$-th root of unity
(i.e., $\eta^k\neq 1$ for $k=0,\ldots,n-1$, and $\eta^n=1$).
We let $\eta=\omega:=\e^\frac{2\pi\ii}{n}$ and $\cF=\cF_\omega$.

This allows us to factorize any circulant matrix in the form
\begin{equation}
C = \cF \Delta \cF^*,
\label{eq:circfac}
\end{equation}
where
$$
\Delta = \diag\left( \hat{C}(1), \hat{C}(\omega), \ldots, \hat{C}(\omega^{n-1})
\right),
$$
and
$$
\hat{C}(\zeta) = \sum_{k=0}^{n-1} c_k \zeta^{-k}
$$
is the discrete Fourier transform of the first column of $C$.
Given the definition of discrete Fourier transform adopted in the $\FFT$
command of Matlab, we have
$$
\deltab := \diag(\Delta) = \FFT(\bm{c}),
$$
being $\bm{c}=(c_0,c_1,\ldots,c_{n-1})^T$ the first column of $C$.

In \smt, a variable \ttt{C} of class \smcirc\ is a record composed by 4 fields.
The field \ttt{C.type} is set to the string `\emph{circulant}', and is a
reminder, present in all \smt\ data types, denoting the kind of the structured
matrix. The first column of the circulant matrix $C$ gives complete information about it, and
is stored in \ttt{C.c}, while \ttt{C.dim} is the dimension $n$.
The field \ttt{C.ev} contains the vector $\deltab$ of the eigenvalues of $C$;
it is computed when the matrix is created and updated every time it is
modified.
This means that the initial allocation of a circulant matrix, as well as some
operations involving it, takes $O(n\log n)$ floating point
operations (\emph{flops}).

For example, an object of class \smcirc\ can be created specifying
its first column, with the command
\begin{quote}
\footnotesize
\begin{verbatim}
C=smcirc([1;2;3;4])
\end{verbatim}
\end{quote}
and it is visualized either as a matrix
\begin{quote}
\footnotesize
\begin{verbatim}
C =
     1     4     3     2
     2     1     4     3
     3     2     1     4
     4     3     2     1
\end{verbatim}
\end{quote}
or showing its record structure
\begin{quote}
\footnotesize
\begin{verbatim}
C =
smcirc object with fields:
    type: 'circulant'
       c: [4x1 double]
     dim: 4
      ev: [4x1 double]
\end{verbatim}
\end{quote}
depending on how the configuration parameter \ttt{display}\ is set; see the
function \ttt{smtconfig} in Section~\ref{sec:general}.
The structure of the object can also be inspected with the command \ttt{get(C)},
independently on the configuration of the package.
If the column vector passed to \smcirc\ is of class \emph{sparse}, this memory
storage class will be preserved in \ttt{C.c}, but not in \ttt{C.ev}.

All operations between circulant matrices have been
implemented, when possible, by
\emph{fast} algorithms, meaning that they require a complexity smaller than the
corresponding unstructured matrix operations.
For example, when the user computes the sum of two circulant matrices with the
command \ttt{A=C+D}, the function \ttt{plus.m} is automatically called, in
order to sum the \ttt{.c} fields, and to update the \ttt{.ev} field of the
resulting object, as follows
$$
\ttt{A.c = C.c + D.c}, \qquad \ttt{A.ev = C.ev + D.ev}.
$$
To multiply a circulant matrix $C$ times a vector $\bm{x}$, we can exploit the
factorization \eqref{eq:circfac} to obtain
\begin{equation}
C\bm{x} = \IFFT\left( \deltab \circ \FFT(\bm{x}) \right),
\label{eq:ctimesx}
\end{equation}
where $\bm{u}\circ\bm{v}=(u_1v_1,\ldots,u_nv_n)^T$ denotes the Schur product of
two vectors.
This requires only 2 \ttt{fft}'s, since the vector $\deltab$ is stored in the
\ttt{.ev} field of the corresponding object \ttt{C}.
In a similar way, the first column of the product of two circulant matrices is
evaluated by the inverse discrete Fourier transform of the product of the
eigenvalues of the two factors.
In all cases, the computation is optimized in terms of complexity.
After performing many operations which update the eigenvalues of an \smcirc\
object, it may be advisable to recompute \ttt{C.ev}, to improve its accuracy; 
if required, this can be done by
\begin{quote}
\footnotesize
\begin{verbatim}
C=smcirc(C.c),
\end{verbatim}
\end{quote}
as the user is not allowed to directly modify the fields of an object belonging
to an \smt\ class.

\begin{table}[htb]
\footnotesize
\centering
\begin{tabular}{l@{~~~~~~~~}c|l@{~~~~~~~~}c}
\hline
\multicolumn{4}{c}{\bf Operators and special characters} \\
\hline
\ttt{plus} & \ttt{A+B}             & \ttt{power} & \ttt{A.\symbol{94}2} \\
\ttt{uplus} & \ttt{+A}             & \ttt{mldivide} & \ttt{A\symbol{92}B} \\
\ttt{minus} & \ttt{A-B}            & \ttt{mrdivide} & \ttt{A/B} \\
\ttt{uminus} & \ttt{-A}            & \ttt{ldivide} & \ttt{A.\symbol{92}B} \\
\ttt{mtimes} & \ttt{A*B}           & \ttt{rdivide} & \ttt{A./B} \\
\ttt{times} & \ttt{A.*B}           & \ttt{transpose} & \ttt{A.'} \\
\ttt{mpower} & \ttt{A\symbol{94}2} & \ttt{ctranspose} & \ttt{A'} \\
\hline
\end{tabular}
\caption{Overloaded operators}
\label{tab:operators}
\end{table}

All the \emph{overloaded} operators, or \emph{methods}, for \smcirc\ objects
are coded in a set of functions, whose names (fixed by the Matlab syntax) are
reported in Table~\ref{tab:operators}, together with the equivalent Matlab
notations.
Each of these functions is called when at least one of the operands in an
expression is of class \smcirc; if the two operands are different \ttt{smt}\
objects, the method corresponding to the first one is called.
The result is structured whenever this is possible.

When an operation is performed between two circulant matrices, the complexity is not larger than
$O(n\log n)$ (for example in the matrix product), while it may be larger when
one of the arguments is unstructured; e.g., the product between a circulant and
a full matrix, which is computed by multiplying the first operand times each of
the columns of the second one, takes $O(n^2\log n)$ \emph{flops}.
So, if \ttt{C} and \ttt{D} are both \smcirc\ objects and \ttt{x} is a vector,
\begin{quote}
\footnotesize
\begin{tabbing}
XXXXXXXXXXXXXXX \= \kill
\ttt{E=C*D} \> produces a \smcirc\ object, \\
\ttt{y=(C+D)*x} \> returns a vector, \\
\ttt{F=C*rand(n)} \> returns an unstructured matrix,
\end{tabbing}
\end{quote}
and, in all cases, the fast algorithms implemented for the \smcirc\ class
are automatically used in the computation.

The operations on \smcirc\ objects which rely on the factorization
\eqref{eq:circfac}, and so exhibit a $O(n^\alpha\log n)$ complexity, are the
product, the power and the left/right division for matrices.
Obviously, some operators do not involve floating point computations at all,
like the transposition or the unary minus.

\begin{figure}[!ht]
\centering
\begin{minipage}{.7\textwidth}
\begin{Algorithm}[The \ttt{plus.m} function for \smcirc\
objects]\label{algo:plus}~\rm
\begin{tabbing}
xxx\=xxxx\=xxxx\=xxxx\=xxxxxxxxxxxxxxxx\=\kill
check validity and dimensions of input arguments $C$ and $D$ \\
deal with scalar or empty arguments \\
\IF $C$ is \texttt{smcirc} \\
\> \IF $D$ is scalar or \smcirc \\
\> \> the result is \smcirc \\
\> \ELSE \IF $D$ is \smtoep \\
\> \> the result is \smtoep \\
\> \ELSE \\
\> \> the result is full \\
\ELSE \` \emph{\{if $C$ is not \emph{\ttt{smcirc}},
	then $D$ is\}} \\
\> \IF $C$ is scalar \\
\> \> the result is \smcirc \\
\> \ELSE \\
\> \> the result is full
\end{tabbing}
\end{Algorithm}
\end{minipage}
\end{figure}

A trivial implementation of the algorithms is not sufficient to obtain
a package which is both robust and transparent for the user.
In fact, each function should be able to handle most of possible user's errors,
and should replicate the typical behaviour of Matlab when any of the operands
are scalars or empty arrays.
As an example, we report in Algorithm~\ref{algo:plus} the structure of the
\ttt{plus.m} function, which is called when the first structured argument in a
sum is an \smcirc\ object.
As it can be seen, when the operands are of different classes, the result
belongs to the less structured class; e.g., circulant plus full is full,
circulant plus Toeplitz is Toeplitz, etc.

Many Matlab standard functions have been redefined for circulant matrices.
They are listed and briefly described in Table~\ref{tab:functions}; among them,
there are simple manipulation and conversion functions, like \ttt{abs},
\ttt{double} or \ttt{full}, some which return logical values (the \ttt{isXXX}
functions), and a few which optimize some computations for \smcirc\ objects,
like \ttt{det}, \ttt{eig} or \ttt{inv}.
The implementation of the last three functions is straightforward, as each
\smcirc\ object contains the eigenvalue of the circulant matrix in its
\ttt{.ev} field.
We remark that some functions require a larger complexity for a
circulant than for a full matrix, like \ttt{imag}, because extracting the
imaginary part of the entries of a circulant matrix requires to recompute its
eigenvalues.

\begin{table}[htb]
\footnotesize
\centering
\begin{tabular}{ll|ll}
\hline
\multicolumn{4}{c}{\bf Elementary math functions} \\
\hline
\ttt{abs}   & absolute value     & \ttt{fix}   & round towards zero\\
\ttt{angle} & phase angle        & \ttt{floor} & round towards $-\infty$\\
\ttt{conj}  & complex conjugate  & \ttt{ceil}  & round towards $+\infty$\\
\ttt{imag}  & imaginary part     & \ttt{round} & round argument\\
\ttt{real}  & real part          & \ttt{sign}  & signum function\\
\hline
\multicolumn{4}{c}{\bf Basic array information} \\
\hline
\ttt{size}    & size of array      &\ttt{get}     & get object fields \\
\ttt{length}  & length of array    &\ttt{isempty}  & true for empty array\\
\ttt{display} & display array      &\ttt{isequal}  & true for equal arrays\\
\hline
\multicolumn{4}{c}{\bf Array operations and manipulation} \\
\hline
\ttt{diag} & diagonals of a matrix  &\ttt{reshape} & change size\\
\ttt{full} & convert to full matrix &\ttt{tril}    & lower triangular part\\
\ttt{prod} & product of elements    &\ttt{triu}    & upper triangular part\\
\ttt{sum}  & sum of elements\\
\hline
\multicolumn{4}{c}{\bf Array utility functions} \\
\hline
\ttt{double}  & convert to double            & \ttt{subsasgn}  & subscripted
assignment\\
\ttt{single}  & convert to single            & \ttt{subsindex} & subscript
index\\
\ttt{isa}     & true if object is in a class & \ttt{subsref}   & subscripted
reference\\
\ttt{isfloat} & true for floating point      & \ttt{end}       & last index\\
\ttt{isreal}  & true for real array          \\ 
\hline
\multicolumn{4}{c}{\bf Matrices and numerical linear algebra} \\
\hline
\ttt{det} & determinant & \ttt{inv} & matrix inverse \\
\ttt{eig} & eigenvalues and eigenvectors \\
\hline
\end{tabular}
\caption{Overloaded functions}
\label{tab:functions}
\end{table}

The list in Table~\ref{tab:functions} is surely incomplete, since in
principle all Matlab matrix functions could be overloaded for circulant
matrices.
We implemented those functions which we consider useful, leaving an extension
of this list, if motivated by real need, to future versions of the package.
It is sufficiently easy to add new methods to the class, since the user can
start from an existing function, as a template, and then place the new file in
the \ttt{smt/@smcirc} directory.

Let us add some comments on some of the functions listed in Table~\ref{tab:functions}.
When adding a new class to Matlab, there are a number of functions which must
be defined so that the class conforms to Matlab syntax rules.
The \ttt{get} method allows to extract a field from an object, while
\ttt{display} defines how an object should be visualized on the screen; this
can be customized in \smt, as it will be shown in Section~\ref{sec:general}.
Some other functions define the effect of subindexing on the new class.
We let two of them, \ttt{subsasgn} and \ttt{subsindex}, just return an error
code for an \smcirc\ object, since we consider them useless for circulant
matrices.
The third one, \ttt{subsref}, is a function which allows to access a field
(\ttt{C.c}) or an element (\ttt{C(2,3)}) of a circulant matrix, and to use
typical Matlab subindexing expressions like \ttt{C(:)} or \ttt{C(3:4,:)}.
Notice that \ttt{C(1:3,4:7)} returns a Toeplitz matrix (i.e., an \smtoep\
object; see Section~\ref{sec:smtoep}), while \ttt{C([1,3,5],6:8)} returns a
full matrix.

The class \smcirc\ includes two additional methods: \ttt{smtvalid} is a
function, called by other functions of the toolbox, which determines if an
object is a valid operand in an expression, while \smtoep\ converts an \smcirc\
object into an \smtoep\ one, as a circulant matrix is also a Toeplitz matrix.

\subsection{Toeplitz matrices}\label{sec:smtoep}

A Toeplitz matrix of order $n$ is a matrix $T$ whose elements are constant
along diagonals, that is
$$
T_{ij} = t_{i-j}, \quad i,j=1,\ldots,n,
$$
e.g., for $n=4$,
\begin{equation}
T = \begin{bmatrix}
t_0 & t_{-1} & t_{-2}& t_{-3}\\
t_1 & t_0 & t_{-1} & t_{-2}\\
t_2 & t_1 & t_0 & t_{-1} \\
t_3 & t_2 & t_1 & t_0
\end{bmatrix}.
\label{eq:toepmat}
\end{equation}
We introduced a class \smtoep, for Toeplitz matrices, similar to the \smcirc\
class.
An \smtoep\ object can be created by specifying its first column and row, for
example with the command
\begin{quote}
\footnotesize
\begin{verbatim}
T=smtoep([4:7],[4:-1:1]),
\end{verbatim}
\end{quote}
or giving only the first column, in which case the resulting matrix is
Hermitian.
Similarly to what happens to \smcirc\ objects, an \smtoep\ object can be
displayed either as
\begin{quote}
\footnotesize
\begin{verbatim}
T =

     4     3     2     1
     5     4     3     2
     6     5     4     3
     7     6     5     4
\end{verbatim}
\end{quote}
or
\begin{quote}
\footnotesize
\begin{verbatim}
T =
smtoep object with fields:
    type: 'toeplitz'
       t: [7x1 double]
    dim1: 4
    dim2: 4
     cev: [8x1 double]
\end{verbatim}
\end{quote}
depending on the \ttt{display}\ configuration parameter; see
\ttt{smtconfig} in Section~\ref{sec:general}.

An \smtoep\ object has two fields for the number of rows (\ttt{T.dim1}) and
columns (\ttt{T.dim2}) of the matrix, while \ttt{T.t} contains the data to
reconstruct the matrix, namely the first row and column, in a form which is
convenient for computation; in the above example,
$$
\ttt{T.t} = (1,2,\ldots,7)^T.
$$
The meaning of the \ttt{T.cev} field will be explained later in this section.

The componentwise operators, like sum, subtraction, and the so-called
dot-operators of Matlab, can be easily implemented for Toeplitz matrices,
similarly to what has been done for the \smcirc\ class.

Regarding the matrix product, it is well known that a Toeplitz matrix can be
embedded in a circulant matrix $C_T$; e.g., given the matrix
\eqref{eq:toepmat}, we can write
\begin{equation}
C_T = \left[ \begin{array}{cccc|cccc}
t_0 & t_{-1} & t_{-2} & t_{-3} & 0 & t_3 & t_2 & t_1 \\
t_1 & t_0 & t_{-1} & t_{-2} & t_{-3} & 0 & t_3 & t_2 \\
t_2 & t_1 & t_0 & t_{-1} & t_{-2} & t_{-3} & 0 & t_3 \\
t_3 & t_2 & t_1 & t_0 & t_{-1} & t_{-2} & t_{-3} & 0 \\
\hline
0 & t_3 & t_2 & t_1 & t_0 & t_{-1} & t_{-2} & t_{-3} \\
t_{-3} & 0 & t_3 & t_2 & t_1 & t_0 & t_{-1} & t_{-2} \\
t_{-2} & t_{-3} & 0 & t_3 & t_2 & t_1 & t_0 & t_{-1} \\
t_{-1} & t_{-2} & t_{-3} & 0 & t_3 & t_2 & t_1 & t_0
\end{array} \right].
\label{eq:ctmat}
\end{equation}
So, to compute the product $\bm{y}=T\bm{x}$ by a fast algorithm, one can
construct a vector $\overline{\bm{x}}$ by padding $\bm{x}$ with zeros to reach
the dimension of $C_T$, then compute $C_T\overline{\bm{x}}$ by
\eqref{eq:ctimesx}, and finally extract $\bm{y}$ from the first components of
the result.

The zero diagonals in \eqref{eq:ctmat} can be deleted, in which case the
dimension of $C_T$ is minimal, or ``tight'': if $T$ is $m\times n$, then the
``tight'' dimension of $C_T$ is $m+n-1$.
On the contrary, we can insert as many zero diagonals as we want.
This may be useful, because the implementations of the \ttt{fft} perform better
when the length of the input vector is a power of 2.

In \smt\ both choices are available, and can be selected by editing the command
\ttt{smtconst}; see Section~\ref{sec:general}.
Although Matlab implementation of the \ttt{fft}, namely FFTW \cite{fftw05},
exhibits a very good performance also when the size of the input vector is a
prime number, we observed that matrix product is generally faster if we extend
the matrix $C_T$ to the next power of 2 exceeding $m+n-1$.

A particular function has been created to speed-up Toeplitz matrix
multiplication.
Thus, the command
\begin{quote}
\footnotesize
\begin{verbatim}
T=toeprem(T)
\end{verbatim}
\end{quote}
pre-computes the eigenvalues of the matrix $C_T$, and stores them in the
\ttt{.cev} field.
This is done automatically when an \smtoep\ object is allocated, and allows to
perform only two \ttt{fft}'s for each matrix product, instead of three.
The price to pay is that, like in the case of circulant matrices, some
elementary functions involving \smtoep\ objects have a complexity larger than
expected, as they need to compute the \ttt{.cev} field of the result.
If this behaviour is not convenient, the automatic call to \ttt{toeprem} can be
disabled by the \ttt{smtconfig} command (see Section~\ref{sec:general}), and the user can either call
\ttt{toeprem} when needed, or renounce to multiplication speedup.

All the operators and functions of Tables~\ref{tab:operators}
and~\ref{tab:functions} have been implemented for the class \smtoep, with some
differences.
Unlike the circulant matrices, there is not a standard method to invert a
Toeplitz matrix, or to compute its determinant or eigenvalues.
On the contrary, various different algorithms are available and, probably, more
will be developed in the future.
For this reason the functions \ttt{inv}, \ttt{det} and \ttt{eig}, supplied with
the toolbox, return an error for an \smtoep\ object, and they are intended to
be overwritten by user supplied programs.

\subsection{Linear systems and preconditioners}\label{sec:Toesys}

Solving circulant linear system is immediate, by employing the factorization
\eqref{eq:circfac}, and the computation requires just two \ttt{fft}'s.
The algorithm is implemented in the functions \ttt{mldivide} and
\ttt{mrdivide}, placed into the \ttt{@smcirc} directory, and is accessible via
the usual matrix left/right division operators.

To solve Toeplitz linear systems, one possibility is to use an iterative solver, either
user supplied, or among those (\ttt{pcg}, \ttt{gmres}, etc.)\ available in
Matlab.
This can be done transparently, taking advantage of the compact storage and
fast matrix-vector product provided by the toolbox.

It is usual to employ preconditioners to speed up the convergence of iterative
methods.
In the case of Toeplitz linear systems, it has been proved that various classes
of circular preconditioners guarantee superlinear convergence for the conjugate
gradient method; see \cite{cj07}.

The function \ttt{smtcprec}, included in the toolbox, provides the
three best known circulant preconditioners, and can be easily extended to
include more.
For a given Toeplitz matrix $A$, the function can construct the Strang preconditioner
\cite{str86}, which suitably modifies the matrix to make it circulant,
the so-called optimal preconditioner \cite{cha88}, which is the solution of
the optimization problem
$$
\min_{C\in\mathcal{C}} \|C-A\|_F,
$$
where $\mathcal{C}$ is the algebra of circulant matrices and $\|\cdot\|_F$
denotes the Frobenius norm, and the superoptimal preconditioner
\cite{cjy91b,tis91,tyr92}, which minimizes
$$
\|I-C^{-1}A\|_F
$$
for $C\in\mathcal{C}$.
While the Strang preconditioner is defined only when $A$ is Toeplitz, the
optimal and superoptimal preconditioners can be computed for any matrix; the
function \ttt{smtcprec} allows this, though the computation of the
preconditioner is fast only for a Toeplitz matrix.

The code in the functions \ttt{strang}, \ttt{optimal}, and \ttt{superopt} (see
Table~\ref{tab:compgen}) was developed by one of the authors during the
research which led to \cite{mrs06a}, and the details of the algorithms are
described in that paper.
To compute the superoptimal preconditioner, it is possible to use either the
method introduced in \cite{cjy91b}, or the one from \cite{tyr92}.
It is remarkable that, since the second algorithm is based on the use of
certain Toeplitz matrices, our implementation is greatly simplified, as it
performs the computation using the arithmetics provided by the toolbox itself.

Using circulant preconditioners with the iterative methods available in Matlab
is straightforward, as these functions use the matrix left division to apply a
preconditioner, and so take advantage of the storage and fast algorithms
furnished by our toolbox.
For example, the instructions
\begin{quote}
\footnotesize
\begin{verbatim}
T=smtgallery('gaussian',5000);
b=T*ones(5000,1);
C=smtcprec('strang',T);
[x,flag,relres,iter]=pcg(T,b,[],[],C);
\end{verbatim}
\end{quote}
create a Gaussian linear system of dimension 5000 (see
Section~\ref{sec:general}) with prescribed solution, and
solve it by the conjugate gradient method, preconditioned by the Strang
circulant preconditioner.

\begin{table}[htb]
\footnotesize
\centering
\noindent\begin{tabular}{ll|ll}
\hline
\multicolumn{2}{c|}{\bf Preconditioners} & \multicolumn{2}{c}{\bf Direct solvers} \\
\hline
\ttt{smtcprec} & circulant preconditioners & \ttt{toms729} & Toeplitz solver \\
\ttt{strang} & Strang preconditioner & \ttt{tlls} & Toeplitz LS solver \\
\ttt{optimal} & optimal preconditioner & \ttt{tsolve} & user supplied function \\
\ttt{superopt} & superoptimal precond. & \ttt{tsolvels} & user supplied function \\
\hline
\multicolumn{4}{c}{\bf General functions} \\
\hline
\ttt{issmcirc} & true for \smcirc\ object   & \ttt{smtconfig} & toolbox configuration  \\
\ttt{issmtoep} & true for \smtoep\ object  & \ttt{smtconst} & toolbox constants setting \\
\ttt{smtcheck} & check toolbox installation  & \ttt{smtgallery} & test matrices \\
\hline
\end{tabular}
\caption{Computational and general functions}
\label{tab:compgen}
\end{table}

Besides the iterative methods, there are also many fast and superfast direct
solvers for a Toeplitz linear system, and some of them have been implemented in
publicly available subroutines.
With our toolbox, we distribute two of them, having computational complexity
$O(\alpha n^2)$; they are called when one of the matrix division operators
(either \ttt{\symbol{92}} or \ttt{/}) are used to invert a Toeplitz matrix.
The related files, listed in Table \ref{tab:compgen},
are placed in the \ttt{private}\ subdirectory of \ttt{smt/@smtoep}.

The first one, \ttt{toms729}, is an implementation of the extended Levinson
algorithm for nonsingular Toeplitz linear systems~\cite{hc92}, written in
Fortran, for which a Matlab MEX gateway is available~\cite{ar08}. This solver,
which has been implemented only for real matrices, calls the \ttt{dsytep}
subroutine from~\cite{hc92} if the system matrix is symmetric, and \ttt{dgetep}
in the general case, with the \ttt{pmax} parameter set to 10.

When the linear system is overdetermined (and full-rank) the toolbox calls the
C-MEX program \ttt{tlls}, developed in~\cite{rod06}, which converts it into a
Cauchy-like system $C\bm{y}=\bm{f}$, and computes its least-squares solution as
the Schur complement of the augmented matrix
$$
M_C=\left[\begin{array}{cc|c}
I & C & 0 \\ C^* & 0 & C^*\bm{f} \\ \hline 0 & I & 0
\end{array}\right],
$$
using the generalized Schur algorithm with partial pivoting.
Complex linear systems are supported.
If the matrix is either underdetermined or rank-deficient, an error is
returned.

It is possible for the user to use different algorithms, by supplying the
functions \ttt{tsolve}, for a nonsingular linear system, or \ttt{tsolvels}, for
least-squares, overwriting those placed in the directory
\ttt{smt/@smtoep/private}, and changing the default behaviour of the toolbox
with the \ttt{smtconfig} command.
For example, entering from the command line the instructions
\begin{quote}
\footnotesize
\begin{verbatim}
smtconfig intsolve off
x=T\b
\end{verbatim}
\end{quote}
disables the solver \ttt{toms729}, and solves the linear system
$T\bm{x}=\bm{b}$ by the user supplied function \ttt{tsolve}.

\subsection{Other functions}\label{sec:general}

Besides the overloaded operators and methods located in the \ttt{@smcirc} and
\ttt{@smtoep} directories, some general functions, listed in
Table~\ref{tab:compgen}, are placed in the main toolbox directory, and are
directly accessible to the user.

Among these functions, we find the two \ttt{isXXX} functions, which return
logical values and check if the supplied parameter belongs to the \ttt{XXX}
class, the \ttt{smtcheck} function, which verifies if the toolbox is correctly
installed, and the function \ttt{smtconst}, intended to define global constants
(for the moment only the dimension of the circulant embedding
\eqref{eq:ctmat}).

Of particular relevance is the \ttt{smtconfig} function, which modifies the
behaviour of the toolbox for what concerns the display method for objects
(\ttt{display} parameter), the use of the Toeplitz premultiplication routine,
discussed in Section \ref{sec:smtoep} (\ttt{toeprem} parameter), the warnings
setting (\ttt{warnings} parameter), and the active Toeplitz solvers
(\ttt{intsolve}\ and \ttt{intsolvels} parameters, see
Section~\ref{sec:Toesys}).
For example,
\begin{quote}
\footnotesize
\begin{tabbing}
XXXXXXXXXXXXXXXXXXXXXXXXX \= \kill
\ttt{smtconfig display compact} (or \ttt{off}) \> selects compact display of objects, \\
\ttt{smtconfig display full} (or \ttt{on}) \> restores standard display method.
\end{tabbing}
\end{quote}
All parameter are set by default to the \ttt{on} state; calling
the command \ttt{smtconfig}\ with no parameters shows the state of all
settings.

\begin{table}[htb]
\footnotesize
\centering
\begin{tabular}{ll}
\hline
\multicolumn{2}{c}{\bf Circulant matrices} \\
\hline
\ttt{crrand}   & uniformly distributed random matrix \\
\ttt{crrandn}  & normally distributed random matrix \\
\hline
\multicolumn{2}{c}{\bf Toeplitz matrices} \\
\hline
\ttt{algdec}    & matrix with algebraic decay \\
\ttt{expdec}    & matrix with exponential decay \\
\ttt{gaussian}  & Gaussian matrix \\
\ttt{tchow}     & Chow matrix \\
\ttt{tdramadah} & matrix of 0/1 with large determinant or inverse \\
\ttt{tgrcar}    & Grcar matrix \\
\ttt{tkms}      & Kac-Murdock-Szego matrix \\
\ttt{tparter}   & Parter matrix \\
\ttt{tphans}    & rank deficient matrix from \cite{han98} \\
\ttt{tprand}    & uniformly distributed random matrix \\
\ttt{tprandn}   & normally distributed random matrix \\
\ttt{tprolate}  & prolate matrix \\
\ttt{ttoeppd}   & symmetric positive definite Toeplitz matrix \\
\ttt{ttoeppen}  & pentadiagonal Toeplitz matrix \\
\ttt{ttridiag}  & tridiagonal Toeplitz matrix \\
\ttt{ttriw}     & upper triangular matrix discussed by Wilkinson \\
\hline
\end{tabular}
\caption{Test matrices available in the \ttt{smtgallery} function}
\label{tab:smtgallery}
\end{table}

We end this section reporting another important feature of the toolbox.
A collection of test matrices is available in \ttt{smtgallery}, which is
modelled on Matlab \ttt{gallery} function, but returns structured objects.
The collection, listed in Table~\ref{tab:smtgallery}, includes random matrices,
three matrices studied in \cite{ms05} (\ttt{algdec}, \ttt{expdec} and
\ttt{gaussian}), one from \cite{han98}, and all the Toeplitz matrices provided
by \ttt{gallery}, most of which come from \cite{hig95}.
The syntax of \ttt{smtgallery} is in the same style as the
Matlab \ttt{gallery} function, and documentation is provided for each test
matrix. For instance,
\begin{quote}
\footnotesize
\begin{tabbing}
XXXXXXXXXXXXXXXXXXXX \= \kill
\ttt{T=smtgallery('gaussian',7)} \> constructs a $7\times 7$ Toeplitz Gaussian
    matrix, \\
\ttt{C=smtgallery('crrand',7,'c')} \> returns a random complex circulant
    matrix, \\
\ttt{help private/tchow} \> displays the help page for the Chow matrix.
\end{tabbing}
\end{quote}

\section{Implementation issues}\label{sec:implementation}

The toolbox is entirely written in the Matlab programming language.
The version officially supported is 7.7 (i.e., release 2008b), anyway we tested
it on previous versions, way back to 6.5, without problems.
To install it, it is sufficient to uncompress the archive file containing the
software, creating in this way the directory \smt\ and its subtree.
This directory must be added to the Matlab search path by the command
\ttt{addpath}, in order to be able to use the toolbox from any other directory.

As noted in the previous sections, a great effort has been devoted to catch all
possible user's errors, and to reproduce the standard behaviour of Matlab,
for example for what concerns the output of each function in the presence of
empty or scalar arrays in input.
Since these features are scarcely documented in Matlab manuals, our choices are
mostly due to experimental tests.

We remark that the \ttt{smtconfig} function, described in
Section~\ref{sec:general}, relies on the use of \emph{warnings}, so issuing the
Matlab command \ttt{warning on} restores the initial configuration of the
toolbox, while \ttt{warning off} may cause unpredictable results.

Some of the toolbox functions use the \ttt{isfloat} command, which was
introduced in version 7 of Matlab.
For those who are using version 6.5, a patch for this function is included in
the software; see the \ttt{README.txt} file in the main toolbox directory.

The two programs used to solve Toeplitz linear system are the only ones which
need to be compiled: \ttt{toms729}~\cite{hc92} was written in Fortran, and uses
the Fortran-MEX gateway from~\cite{ar08}, while \ttt{tlls}~\cite{rod06} was
originally developed as a C-MEX program.
The MEX interface is a library, distributed with Matlab, which allows to
dynamically link a Fortran or C subroutine to Matlab, and to exchange input and
output parameters between the compiled program and the environment, using the
usual Matlab syntax.

Both Toeplitz solvers can be easily compiled under Linux, using the
\ttt{Makefile} placed in the directory \ttt{smt/@smtoep/private}, and we
provide precompiled executables for 32 and 64 bits architectures.
Compiling the same programs under Windows is a bit more involved: we used a
porting of the GNU-C compiler~\cite{mingw} and the ``MEX configurator''
Gnumex~\cite{gnumex}, but precompiled executables are available for various
Matlab versions; see the \ttt{README.txt} file and the content of the
\ttt{smt/@smtoep/private} directory.

\section{Conclusion}\label{sec:conclusion}

The Structured Matrix Toolbox is a Matlab package which implements optimized
storage and fast arithmetics for circulant and Toeplitz matrices, offering a
robust and easily extensible framework.
The toolbox is available at the web page \url{http://bugs.unica.it/smt/}.
We are currently performing numerical tests to assess its performance, and
there are plans to extend its functionality by adding the support for other
classes of structured matrices.


\begin{thebibliography}{10}

\bibitem{ar08}
A.~Aric\`o and G.~Rodriguez.
\newblock {\em \texttt{toms729gw}: a Matlab (Fortran) MEX Gateway for TOMS
  Algorithm 729, by P. C Hansen and T. Chan}.
\newblock University of Cagliari, 2008.
\newblock Available at: \url{http://bugs.unica.it/~gppe/soft/}.

\bibitem{cj07}
R.~H. Chan and X.-Q. Jin.
\newblock {\em An Introduction to Iterative Toeplitz Solvers}.
\newblock SIAM, Philadelphia, 2007.

\bibitem{cjy91b}
R.~H. Chan, X.-Q. Jin, and M.~C. Yeung.
\newblock The circulant operator in the {B}anach algebra of matrices.
\newblock {\em Linear Algebra Appl.}, 149:41--53, 1991.

\bibitem{cha88}
T.~F. Chan.
\newblock An optimal circulant preconditioner for {T}oeplitz systems.
\newblock {\em SIAM J. Sci. Stat. Comput.}, 9(4):766--771, 1988.

\bibitem{dav79}
P.~J. Davis.
\newblock {\em Circulant Matrices}.
\newblock Wiley, New York, 1979.

\bibitem{dd05}
J.~Demmel and J.~Dongarra.
\newblock {ST-HEC}: Reliable and scalable software for linear algebra
  computations on high end computers.
\newblock Available at:
  \url{http://www.cs.berkeley.edu/~demmel/Sca-LAPACK-Proposal.pdf}, 2005.

\bibitem{fftw05}
M.~Frigo and S.G. Johnson.
\newblock The design and implementation of {FFTW3}.
\newblock {\em IEEE Proc.}, 93(2):216--231, 2005.
\newblock Software available at: \\\url{http://www.fftw.org/}.

\bibitem{han94}
P.~C. Hansen.
\newblock Regularization tools: A {M}atlab package for analysis and solution of
  discrete ill-posed problems.
\newblock {\em Numer. Algorithms}, 6:1--35, 1994.
\newblock Software available at: \url{http://www2.imm.dtu.dk/~pch/Regutools/}.

\bibitem{han98}
P.~C. Hansen.
\newblock {\em Rank-Deficient and Discrete Ill-Posed Problems: Numerical
  Aspects of Linear Inversion}.
\newblock SIAM Monographs on Mathematical Modeling and Computation. SIAM,
  Philadelphia, 1998.

\bibitem{hc92}
P.~C. Hansen and T.~F. Chan.
\newblock Fortran subroutines for general {T}oeplitz systems.
\newblock {\em ACM Trans. Math. Software}, 18(3):256--273, 1992.

\bibitem{moore}
P.~C. Hansen, M.~Jacobsen, and T.~K. Jensen.
\newblock {\em MOORe Tools: Modular Object Oriented Regularization Tools}.
\newblock Technical University of Denmark, Informatics and Mathematical
  Modelling, 2006.
\newblock Available at: \url{http://www2.imm.dtu.dk/~pch/MOOReTools/}.

\bibitem{hig95}
N.~J. Higham.
\newblock The {T}est {M}atrix {T}oolbox for {M}atlab ({V}ersion 3.0).
\newblock Technical Report 276, Manchester Centre for Computational
  Mathematics, 1995.

\bibitem{matlab77}
The MathWorks, Natick.
\newblock {\em Matlab ver.\ 7.7}, 2008.

\bibitem{restoretools}
J.~Nagy.
\newblock {\em RestoreTools, An Object Oriented Matlab Package for Image
  Restoration}.
\newblock Emory University, Department of Mathematics and Computer Science,
  2007.
\newblock Available at: \\\url{http://www.mathcs.emory.edu/~nagy/RestoreTools/}.

\bibitem{rod06}
G.~Rodriguez.
\newblock Fast solution of {T}oeplitz- and {C}auchy-like least squares
  problems.
\newblock {\em SIAM J. Matrix Anal. Appl.}, 28(3):724--748, 2006.

\bibitem{gnumex}
SourceForge.net.
\newblock {\em Gnumex}, 2000.
\newblock Available at: \\\url{http://gnumex.sourceforge.net/}.

\bibitem{mingw}
SourceForge.net.
\newblock {\em MinGW}, 2008.
\newblock Available at: \\\url{http://www.mingw.org/}.

\bibitem{str86}
G.~Strang.
\newblock A proposal for {T}oeplitz matrix calculations.
\newblock {\em Stud. Appl. Math.}, 74:171--176, 1986.

\bibitem{tis91}
M.~Tismenetsky.
\newblock A decomposition of {T}oeplitz matrices and optimal circulant
  preconditioning.
\newblock {\em Linear Algebra Appl.}, 154/156:105--121, 1991.

\bibitem{tyr92}
E.~E. Tyrtyshnikov.
\newblock Optimal and superoptimal circulant preconditioners.
\newblock {\em SIAM J. Matrix Anal. Appl.}, 13(2):459--473, 1992.

\bibitem{maseteam}
M.~Van~Barel.
\newblock {\em Software produced by members of the MaSe-Team (Matrices having
  Structure)}.
\newblock Katholieke Universiteit Leuven, Department of Computer Science, 2008.
\newblock Available at: \\\url{http://www.cs.kuleuven.ac.be/~marc/software/}.

\bibitem{mrs06a}
C.V.M. {van der}~Mee, G.~Rodriguez, and S.~Seatzu.
\newblock Fast computation of two-level circulant preconditioners.
\newblock {\em Numer. Algorithms}, 41(3):275--295, 2006.

\bibitem{ms05}
C.V.M. {van der}~Mee and S.~Seatzu.
\newblock A method for generating infinite positive self-adjoint test matrices
  and riesz bases.
\newblock {\em SIAM J. Matrix Anal. Appl.}, 26(4):1132--1149, 2005.

\bibitem{slicot}
The Working Group on Software (WGS).
\newblock {\em The Control and Systems Library SLICOT}, 1998.
\newblock Available at: \url{http://www.slicot.org/}.

\bibitem{niconet}
The Working Group on Software (WGS).
\newblock {\em The Numerics in Control Network NICONET}, 1998.
\newblock Available at: \\\url{http://www.icm.tu-bs.de/NICONET/}.

\bibitem{tsmd}
The Working Group on Software (WGS).
\newblock {\em Basic Software Tools for Structured Matrix Decompositions and
  Perturbations}, 2001.
\newblock Available at: \url{http://www.icm.tu-bs.de/NICONET/NICtask1B.html}.

\end{thebibliography}

\end{document}